\parindent=0pt
\magnification=1200

\font\goth=eufm10 
\font\cir=wncyr10

\def\ddefi#1{

{{\bf Definition #1.}

}}
\def\nota#1{

{{\bf Notation #1.}

}}
\def\lem#1{

{{\bf Lemma #1.}

}}
\def\prop#1{

{{\bf Proposition #1.}

}}
\def\teo#1{

{{\bf Theorem #1.}

}}
\def\cor#1{

{{\bf Corollary #1.}

}}
\def\dim{

{{\bf Proof:} }}
\def\sqr#1#2{{\vcenter{\vbox{\hrule height .#2pt
                             \hbox{\vrule width .#2pt height#1pt\kern#1pt
                                   \vrule width .#2pt}
                             \hrule height .#2pt}}}}
\def\square{\mathchoice\sqr54\sqr54\sqr{4.1}3\sqr{3.5}3}
\def\finedim{{\hfill\hbox{\enspace{\bf $\square$}}} 

}
\def\rem#1{

{{\bf Remark #1.}

}}

\def\bbuildrel#1_#2^#3{\mathrel{\mathop{\kern 0pt#1}\limits_{#2}^{#3}}}

\def\rapp#1#2{ { {\displaystyle {#1}} \over {\displaystyle {#2}} }}
\def\ugdef
{{\hskip 2 pt}\dot{=}_{{\hskip -5.5 pt}{\displaystyle .}}{\hskip 5
pt}}

\def\mn{{\bf N}}
\def\mz{{\bf Z}}

\def\mc{{\bf C}}

\def\u{{\cal{U}}}

\def\z{{\cal Z}}
\def\ca{{\cal A}}
\def\hr{{\cal H}^r}
\def\a{\alpha}
\def\d{\delta}
\def\e{\varepsilon}

\def\ug{\underline{\gamma}}

\def\Pp{{\cal P}}

\def\tn{\tilde n}

\def\tPhi{\tilde\Phi}

\def\tk{k}
\def\cb{{\cal B}}

\def\beck{[1]} 
\def\rmcen{[2]} 
\def\rcin{[3]} 
\def\perf{[4]} 
\def\deckac{[5]} 
\def\KM{[6]} 
\def\lusg{[7]} 
\def\lusz{[8]} 
\def\tan{[9]}

\def\intr{0}

\def\taqae{1} 

\def\takman {\taqae.{1}} 

\def\qa {\taqae.{2}} 
\def\uqgr{\qa.{1}} 
\def\qstrutt{\qa.{2}} 

\def\speceps {\taqae.{3}} 
\def\intf{\speceps.{1}} 
\def\altreps{\speceps.{2}} 
\def\lonepr{\speceps.{3}} 
\def\queresta{\speceps.{4}}

\def\ubdz{2} 
\def\ziz{\ubdz.{1}} 
\def\zelem{\ubdz.{2}} 
\def\hypjf{\ubdz.{3}} 
\def\multg{\ubdz.{4}} 
\def\rcoinc{\ubdz.{5}} 
\def\twimmult{\ubdz.{6}} 
\def\checerc{\ubdz.{7}}

\def\zupiu{3} 
\def\qlcnt{\zupiu.{1}} 
\def\adnd{\zupiu.{2}} 
\def\esd{\zupiu.{3}} 
\def\cncls{\zupiu.{4}}

\def\zue{4} 
\def\zk{\zue.{1}} 
\def\zgen{\zue.{2}} 
\def\centrofine{\zue.{3}}

\def\ref{5}

{\centerline{\bf{THE CENTER OF TWISTED AFFINE QUANTUM
ALGEBRAS}}}

{\centerline{\bf{AT ODD ROOTS OF 1}}} 

{\centerline{Ilaria Damiani}}

{\hfill
{
{\it{To the people of Yugoslavia}}}}

{\hfill
{
{\it{for the price they are paying}}}}

{\hfill
{
{\it{to the ``new world order''}}}}

{\hfill
{
{\cir{smrt fashizmu - sloboda narodu}}}}

\input xy
\xyoption{all}

{\bf{\intr. ${\scriptstyle{\rm{\bf INTRODUCTION.}}}$}}

The aim of this paper is describing the center of the
specialization (or, we should say, of a chosen specialization, see
paragraph \speceps) of twisted affine quantum algebras at primitive odd
roots of 1 (with some further very slight conditions on the order of these
roots, see notation
\lonepr). This result is already known for quantum algebras of finite type
(see
\deckac) and of untwisted affine type (see \rmcen); on the other hand the
investigation of the structure of twisted affine quantum algebras carried
out in \rcin\ allows to complete this program with just a little more
effort. 

Actually the structure needed on $\u_q$ is very rich: for this reason 
section \taqae, where the notations are introduced, is the occasion for
recalling, as shortest as possible, the instruments that will be used in
the arguments presented. 

In section \ubdz\ the contravariant form is used, following an argument
introduced in \deckac\ and refined in \rmcen, in order to get a control on
the ``dimension'' of the center: the contravariant form is studied by 
making use of both its general properties and its connections with the
Killing form, which, on its side, was studied in details, for the twisted
algebras, in \rcin. 

In section \zupiu\ some specific computations are developed, in order
to find the missing central elements, which are not simply powers of
root vectors, but which are linear combinations of (non central)
imaginary root vectors. Thanks to these additional elements we obtain the
complete picture of the positive part of $\z(\u_{\e})$. 

Finally section \zue\ is devoted to glueing the pieces together, in order
to state the desired assertion about the center: it is the quotient of an
algebra of polynomials in an infinite number of variables by a very
``small'' ideal, which is indeed a principal ideal generated by an element
$P_Z\in\u_{\e}^0$.

I want to thank Prof. Neda Bokan, who made it possible, as the Dean of the
Faculty of Mathematics of the University of Belgrade, to establish a
scientific cooperation with the Department of Mathematics of the
University of Rome ``Tor Vergata'' and to make this collaboration official
with the signature (by both Rectors, Prof.
Puri\'c and Prof. Finazzi-Agr\`o) of a bilateral agreement. It is thanks to
her efforts that the first contacts between our universities, in summer
1999, immediately after the aggression of the NATO countries against the
Federal Republic of Yugoslavia, turned into proficuous meetings in
Belgrade and into the participation of a delegation from the Second
University of Rome in the X Congress of Yugoslav Mathematicians: for me it
has been a particular pleasure to be in Belgrade again, to meet the
colleagues that I had already met in 1999 and to give a continuity to
our scientific exchanges, which shall go on, in the very next future,
with a program of invitations in Rome. 

As a mathematician concerned about the deformed use of science (to make
war, to destroy a country, to isolate a community,...) and as part of the
movement which in Italy carried out a total opposition against that war
and is now fighting against its continuation in these times of so-called
peace, I'm happy to bring to this Congress my, unfortunately too small,
tribute and solidarity for the price that Yugoslavia paid and is still
paying to the strategies of domination of the world. 

{\bf{\taqae. ${\scriptstyle{\rm{\bf TWISTED\ \ AFFINE\ \ QUANTUM\
\ ALGEBRAS\ \ AT
\ \ ROOTS\ \ OF\ \ 1:\ NOTATIONS.}}}$}}

{\S \takman. {\it The Kac-Moody algebra.}}

Twisted affine quantum algebras are the quantization of the enveloping
algebra of a class of Lie algebras, namely the class of twisted affine
Kac-Moody algebras. 

A complete description of these KM-algebras, as well
as a motivation for their denomination, can be found in \KM, where they
were introduced. What it is important to recall here is that a KM-algebra 
{\goth g} is a Lie algebra whose generators and relations can be expressed
in terms of a matrix $A$, called the (generalized) Cartan matrix of 
{\goth g}, and that the same information contained in $A$ can be encoded
in a diagram $\Gamma$ (the Dynkin diagram of {\goth g}). 

The Dynkin diagrams associated to twisted affine KM-algebras
consist in three families ($A_{2n}^{(2)}$,
$A_{2n-1}^{(2)}$,
$D_{n+1}^{(2)}$) and two isolated cases ($E_{6}^{(2)}$, $D_{4}^{(3)}$),
which are listed below (the type is in general denoted by 
$X_{\tn}^{(k)}$, which means that, for example, in case $E_6^{(2)}$ we
have $X=E$, $\tn=6$, $k=2$):
\vskip -.8 truecm
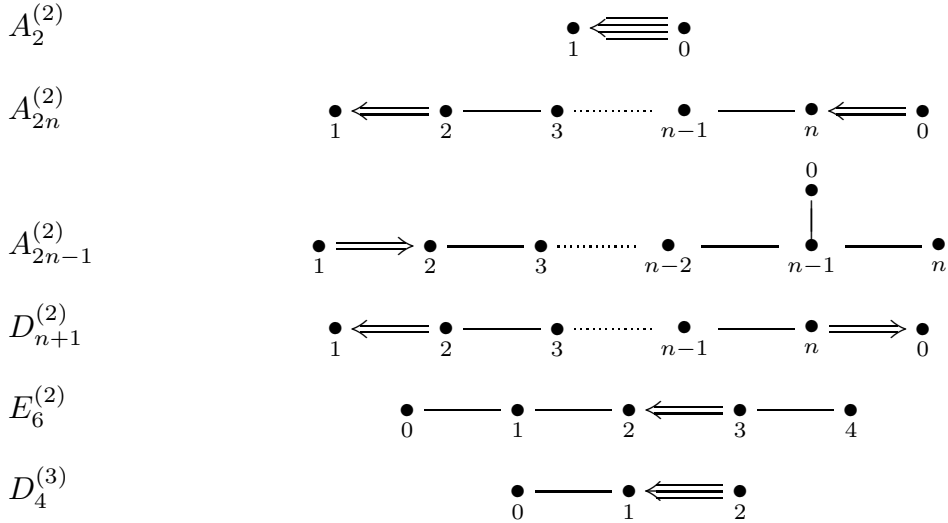
$$
\xymatrix{
{\bbuildrel\bullet_{1}^{{\phantom{1}}}}
&
{\bbuildrel\bullet_{0}^{{\phantom{0}}}}  \ar@{-}@<0.75ex>[l]|>{\quad }
\ar@{=>}[l]
\ar@{-}@<-0.75ex> [l] |>{\quad}
}
\leqno{A_{2}^{(2)}}$$
\vskip -1.0 truecm
$$\xymatrix{{\bbuildrel\bullet_1^{{\phantom{1}}}}&
{\bbuildrel\bullet_2^{{\phantom{2}}}}\ar@{=>}[l]\ar@{-}[r]&
{\bbuildrel\bullet_3^{{\phantom{3}}}}\ar@{.}[r]&
{\bbuildrel\bullet_{n-1}^{{\phantom{n-1}}}}\ar@{-}[r]&
{\bbuildrel\bullet_n^{{\phantom{n}}}}&
{\bbuildrel\bullet_0^{{\phantom{0}}}}\ar@{=>}[l]}\leqno{A_{2n}^{(2)}}$$
\vskip -1.2 truecm
$$\xymatrix{{\phantom{{\bbuildrel\bullet_1^1}}}\ar@{}[r]&
{\phantom{{\bbuildrel\bullet_2^2}}}\ar@{}[r]&
{\phantom{{\bbuildrel\bullet_3^3}}}\ar@{}[r]&
{\phantom{{\bbuildrel\bullet_{n-2}^{n-2}}}}\ar@{}[r]&
\buildrel 0\over\bullet\ar@{}[r]&
{\phantom{{\bbuildrel\bullet_n^n}}}}\leqno{{\phantom{A_{2n-1}^{(2)}}}}$$
\vskip -1.6 truecm
$$\xymatrix{{\phantom{{\bbuildrel\bullet_1^1}}}\ar@{}[r]&
{\phantom{{\bbuildrel\bullet_2^2}}}\ar@{}[r]&
{\phantom{{\bbuildrel\bullet_3^3}}}\ar@{}[r]&
{\phantom{{\bbuildrel\bullet_{n-2}^{n-2}}}}\ar@{}[r]&
|\ar@{}[r]&
{\phantom{{\bbuildrel\bullet_n^n}}}}\leqno{{\phantom{A_{2n-1}^{(2)}}}}$$
\vskip -2.0 truecm
$$\xymatrix{{\phantom{{\bbuildrel\bullet_1^1}}}\ar@{}[r]&
{\phantom{{\bbuildrel\bullet_2^2}}}\ar@{}[r]&
{\phantom{{\bbuildrel\bullet_3^3}}}\ar@{}[r]&
{\phantom{{\bbuildrel\bullet_{n-2}^{n-2}}}}\ar@{}[r]&
|\ar@{}[r]&
{\phantom{{\bbuildrel\bullet_n^n}}}}\leqno{{\phantom{A_{2n-1}^{(2)}}}}$$
\vskip -1.8 truecm
$$\xymatrix{{\bbuildrel\bullet_1^{{\phantom{1}}}}\ar@{=>}[r]&
{\bbuildrel\bullet_2^{{\phantom{2}}}}\ar@{-}[r]&
{\bbuildrel\bullet_3^{{\phantom{3}}}}\ar@{.}[r]&
{\bbuildrel\bullet_{n-2}^{{\phantom{n-2}}}}\ar@{-}[r]&
\bbuildrel\bullet_{n-1}^{{\phantom{n-1}}}\ar@{-}[r]&
{\bbuildrel\bullet_n^{{\phantom{n}}}}}\leqno{A_{2n-1}^{(2)}}$$
\vskip -1.0 truecm
$$\xymatrix{{\bbuildrel\bullet_1^{{\phantom{1}}}}&
{\bbuildrel\bullet_2^{{\phantom{2}}}}\ar@{=>}[l]\ar@{-}[r]&
{\bbuildrel\bullet_3^{{\phantom{3}}}}\ar@{.}[r]&
{\bbuildrel\bullet_{n-1}^{{\phantom{n-1}}}}\ar@{-}[r]&
{\bbuildrel\bullet_n^{{\phantom{n}}}}\ar@{=>}[r]&
{\bbuildrel\bullet_0^{{\phantom{0}}}}}\leqno{D_{n+1}^{(2)}}$$
\vskip -1.0 truecm
$$\xymatrix{{\bbuildrel\bullet_0^{{\phantom{0}}}}\ar@{-}[r]&
{\bbuildrel\bullet_1^{{\phantom{1}}}}\ar@{-}[r]&
{\bbuildrel\bullet_2^{{\phantom{2}}}}&
{\bbuildrel\bullet_3^{{\phantom{3}}}}\ar@{=>}[l]\ar@{-}[r]&
{\bbuildrel\bullet_4^{{\phantom{4}}}}}\leqno{E_{6}^{(2)}}$$
\vskip -1.0 truecm
$$\xymatrix{{\bbuildrel\bullet_0^{{\phantom{0}}}}\ar@{-}[r]&
{\bbuildrel\bullet_1^{{\phantom{1}}}}\ar@3{<-}[r]&
{\bbuildrel\bullet_2^{{\phantom{2}}}}}\leqno{D_{4}^{(3)}}$$
\vskip -.3 truecm

The matrix $A=(a_{ij})_{ij\in I}$ can be recovered from $\Gamma$ as
follows: 

a) $I$ is the set of vertices of $\Gamma$; 

b) $\cases{a_{ii}\!=2&$\forall i\in I$\cr
a_{ij}=-\#\{{\rm{edges\ in\ }}\Gamma\ {\rm{connecting}}\ i\
{\rm{and}}\  j\}&$\forall i\neq j\in I$ s.t. 
$\exists$an arrow pointing at $i$\cr
&or there is no edge connecting $i$ and $j$\cr
a_{ij}\!=-1&otherwise.}$ 

Remark that the 
set $I$ has been identified with $\{0,1,...,n\}$; the set
$\{1,...,n\}\ugdef I\setminus\{0\}$ is denoted by $I_0$. 

Attached to these data there is the notion of root system $\Phi\subseteq
Q\ugdef\oplus_{i\in I}\mz\a_i$, of positive and negative, real and
imaginary roots, and of multiplicity of a root: 
$\Phi_+\ugdef\Phi\cap Q_+=\Phi\cap(\sum_{i\in I}\mn\a_i)$,
$\Phi_-\ugdef-\Phi_+=\Phi\setminus\Phi_+$,
$\Phi^{{\rm{im}}}\ugdef\mz\d\setminus\{0\}$ with 
$\d=\sum_{i\in I}r_i\a_i$, $r_0=1$ and $\sum_{j\in I}a_{ij}r_j=0$ $\forall
i\in I$, and $\Phi^{{\rm{re}}}=\Phi\setminus\Phi^{{\rm{im}}}$. The
multiplicity of each real root is 1, while the multiplicity of $r\d$ can be
described as follows: 
$\forall r\neq 0$ denote by
$I^r$ the set $I^r\ugdef\cases{I_0&in case $A_{2n}^{(2)}$\cr
\{i\in
I_0|d_i|r\}&otherwise}$ where
$\{d_i|i\in I\}$ is the set of positive integers such that $min\{d_i|i\in
I\}=1$ and
$d_ia_{ij}=d_ja_{ji}$
$\forall i,j\in I$ (the $d_i$'s exist and are uniquely determined since
$A$ is symmetrizable and indecomposable, see \KM). Then the multiplicity of
$r\d$ is
$\#I^r$. 

Thus the set of positive roots with multiplicities, denoted by $\tPhi_+$,
can be described as $\tPhi_+=\Phi_+^{{\rm{re}}}\cup\tPhi_+^{{\rm{im}}}=
\Phi_+^{{\rm{re}}}\cup(\cup_{r>0}\{r\d\}\times I^r)$: remark that the
condition $i\in I^r$ is equivalent to $(r\d,i)\in\tPhi_+$ and that 
$\#I^r=\cases{n&if $k|r$ 
\cr
{\tn-n\over k-1}&otherwise.}$

It is worth noticing that $\tPhi_+$ is an index set for a basis of {\goth
n}$_+\subseteq${\goth g} (see \KM).

{\S \qa. {\it The quantum algebra.}}

In this paragraph a short account of the quantum algebra $\u_q$
associated to a (twisted) affine Cartan matrix $A$ (or to the
corresponding Dynkin diagram $\Gamma$) will be given, and some of the main
structures will be shortly recalled.
\ddefi{\uqgr}
$\u_q$ is the $\mc(q)$-associative unitary algebra generated by 
$\{E_i,F_i,K_i^{\pm 1}|i\in I\}$ with relations: 
\vskip -.8 truecm
$$[K_i,K_j]=0,\ K_iE_j=q^{d_ia_{ij}}E_jK_i,\
K_iF_j=q^{-d_ia_{ij}}F_jK_i,\
[E_i,F_j]=\d_{ij}\rapp{K_i-K_i^{-1}}{q^{d_i}-q^{-d_i}}\ \forall i,j\in I$$
\vskip -1 truecm
$$\sum_{r=0}^{1-a_{ij}}\!(-1)^r\!{1-a_{ij}\brack r}_{q^{d_i}}
\!E_i^rE_jE_i^{1-a_{ij}-r}\!=\!0\!=\!\!
\sum_{r=0}^{1-a_{ij}}\!(-1)^r\!{1-a_{ij}\brack r}_{q^{d_i}}\!
F_i^rF_jF_i^{1-a_{ij}-r}\ \forall i\neq j\in I$$
\vskip -.5 truecm
where $\forall m\geq m^{\prime}, r\in\mn$
$[m]_{q^r}\ugdef{q^{mr}-q^{-mr}\over q^r-q^{-r}}$,
$[m]_{q^r}!\ugdef\prod_{s=1}^m[s]_{q^r}$, ${m\brack
m^{\prime}}_{q^{r}}\ugdef
{[m]_{q^r}!\over[m^{\prime}]_{q^r}![m-m^{\prime}]_{q^r}!}$.
In the following remark some fundamental structures that $\u_q$ can be
endowed with are listed, and the references given.
\rem{\qstrutt}
a) $\u_q=\u_q^-\otimes_{\mc(q)}\u_q^0\otimes_{\mc(q)}\u_q^+$ where
$\u_q^-$, $\u_q^0$, $\u_q^+$ are the subalgebras of $\u_q$
respectively generated by 
$\{F_i|i\in I\}$, $\{K_i^{\pm 1}|i\in I\}$, $\{E_i|i\in I\}$ (see \lusg); 

b) $\u_q=\oplus_{\eta\in Q}\u_{q,\eta}$ is the $Q$-gradation induced by
$E_i\in\u_{q,\a_i}$, $F_i\in\u_{q,-\a_i}$, $K_i^{\pm 1}\in\u_{q,0}$; 

c) $\Omega:\u_q\rightarrow\u_q$ is the antilinear antiinvolution such that
$E_i\buildrel\Omega\over\leftrightarrow F_i$,
$K_i\buildrel\Omega\over\leftrightarrow K_i^{-1}$,
$q\buildrel\Omega\over\mapsto q^{-1}$; 

d) the braid group $\cb$ acts on $\u_q$ and its quotient $W$ (the Weyl
group) acts on $Q$ in such a way that $T(\u_{q,\eta})=\u_{q,w(\eta)}$
where $w$ is the image of $T$ in $W$ (see \lusg);

e) $W.\{\a_i|i\in I\}=\Phi^{{\rm{re}}}$; this allows to construct
elements $E_{\a}\in\u_{q,\a}^+$ $\forall\a\in\Phi_+^{{\rm{re}}}$ (with
$E_{\a_i}=E_i$), which are called positive real root vectors (see \lusg,
\beck, \rcin); 

f) by the commutation of the $E_{\a}$'s ($\a\in\Phi_+^{{\rm{re}}}$) the
imaginary root vectors $E_{(r\d,i)}\in\u_{q,r\d}^+$ with
$(r\d,i)\in\tPhi_+^{{\rm{im}}}$ are defined: they generate a commutative
subalgebra of $\u_q^+$ (see \beck, \rcin);

g) given a total ordering $<$ of $\tPhi_+$ set, $\forall\eta\in Q_+$, 
$\Pp(\eta)\ugdef\!\{\ug=(\gamma_1\!\leq\!...\!\leq\!\gamma_m)|
\sum_{s=1}^m\!p(\gamma_s)=\eta\}$ (where $p:\!\tPhi_+\!\rightarrow\!\Phi_+$
is the natural projection),
$par(\eta)\ugdef\!\#\Pp(\eta)$ and, $\forall x\!:\!\tPhi_+\!\ni\a\mapsto
x_{\a}\in\u_q$, 
$\forall\ug=(\gamma_1,...,\gamma_m)\in\Pp(\eta)$, 
$x(\ug)\ugdef x_{\gamma_1}\cdot...\cdot x_{\gamma_m}$; then 
$\{E(\ug)|\ug\in\cup_{\eta\in Q_+}\Pp(\eta)\}$ is a basis of $\u_q^+$,
called PBW-basis (for a discussion on the admitted orderings on
$\tPhi_+$, see \rcin);

h) the Killing form $(\cdot,\cdot):\u_q^{\geq 0}\times\u_q^{\leq
0}\rightarrow\mc(q)$ connects the algebra and coalgebra structures of
$\u_q$ (see \tan); 

i) the (``universal'') contravariant form
$H:\u_q^-\times\u_q^-\rightarrow\u_q^0$ contains information about the
commutation relations (see \deckac); its restriction to $\u_{q,-\eta}^-$
is denoted by $H_{\eta}$.

{\S \speceps. {\it The specialization at (odd) roots of 1.}}

The aim of this
paragraph is to define the specialization of
$\u_q$ at $\e$ when $\e$ is a nonzero complex number. 
Roughly speaking, the idea is to transform the indeterminate
$q$ into a parameter which takes values in $\mc^*$, thus obtaining, for
each value $\e$ of the parameter, an associative $\mc$-algebra, that
will be denoted by $\u_{\e}$; this can be done in different 
ways: in this paper one of these different methods to specialize $\u_q$ is
chosen (see \deckac), but it can be useful to recall that other
specializations (see
\lusz) have an essentially different behaviour exactly under the aspect
that we are going to consider, the center. 
\ddefi{\intf}
Let $\u_q$ be a quantum algebra of type $X_{\tn}^{(k)}$ and let $\ca$ be
the subalgebra of
$\mc(q)$ generated over
$\mc$ by 
$\{q,q^{-1},(q^m-q^{-m})^{-1}|1\leq m\leq k\}$, i.e.
$\ca\ugdef\mc[q,q^{-1},(q^m-q^{-m})^{-1}|1\leq m\leq k]$; the
$\ca$-subalgebra $\u_{\ca}$ of $\u_q$ generated by $\{E_i,F_i,K_i^{\pm
1}|i\in I\}$ is called the integer form of $\u_q$. Remark that if
$\e\in\mc^*$ is such that $\e^{2m}\neq 0$ $\forall m=1,...,k$ than $q-\e$
is not invertible in $\u_{\ca}$. In this case we call specialization of
$\u_q$ at $\e$, and denote it by $\u_{\e}$, the quotient 
$\u_{\e}\ugdef\u_{\ca}/(q-\e)$. 
\rem{\altreps}
It is possible to avoid the restrictions on $\e$ introduced for the
construction above, but since we are actually interested in further
restrictions it is useless, for the purpose of this paper, to look for an
unnecessary generality.
\nota{\lonepr} 
From now on $l$ will denote an odd integer bigger than $k$, and $\e$ a
primitive $l^{{\rm{th}}}$ root of 1. Remark that such an $\e$ satisfies
the restrictions required in definition \intf.
\prop{\queresta}
$\u_{\e}$ inherits from $\u_q$ most of its structures: 

a) the triangular decomposition:
$\u_{\e}=\u_{\e}^-\otimes_{\mc}\u_{\e}^0\otimes_{\mc}\u_{\e}^+$; 

b) the $Q$-gradation: $\u_{\e}=\oplus_{\eta\in Q}\u_{\e,\eta}$; 

c) the antilinear antiinvolution $\Omega:\u_{\e}\rightarrow\u_{\e}$ such
that $\u_{\e}^+{\buildrel\Omega\over\leftrightarrow}\u_{\e}^-$;

d) The structure of $\u_{\e}^0$: $\u_{\e}^0=\mc[K_i^{\pm 1}|i\in I]$;
$\forall\lambda=\sum_{i\in I}m_i\a_i\in Q$, 
$K_{\lambda}\ugdef\prod_{i\in I}K_i^{m_i}$;

e) the braid group action, the root vectors and their basis-properties (PBW
basis).

Moreover recall that $(q-\e)|[m]_{q^r}$ (that is $[m]_{q^r}=0$ in
$\u_{\e}$)
$\Leftrightarrow l|mr$ and
$l\not| r$, and that $(q-\e)^2$ never divides $[m]_{q^r}$ (if $m\neq 0$);
given an element $a$ in an $\ca$-algebra, $mult_{\e}a$ denotes the
multiplicity of $\e$ in $a$ (that is $(q-\e)^{mult_{\e}a}||a$). 
\dim
See \perf, paragraph 1.D1.\finedim

{\bf{\ubdz. ${\scriptstyle{\rm{\bf AN\ \ UPPER\ \ BOUND\ \ FOR\ \ }}}$
$\rm{{\bf dim}{{{(\z(\u_{\e})\cap\u_{\e,\eta}^+).}}}}$}}

The strategy followed to describe the center of the algebras that we are
dealing with passes through the investigation of its positive part,
which is the goal of this and the next section; it consists in exhibiting
a family of central elements, proving that they constitute a set of
generators and showing that there are no relations among them. 

Exhibiting a first set of central elements is not difficult: it can be
done thanks to the results on the commutation relations between the root
vectors and to the commutation properties of the elements $K_{\lambda}$'s. 
Analogously, it is thanks to the PBW basis (see \rcin,
section 6) that we can deduce that there are no relations among these
central root vectors. 

Thus the starting point of this section is an estimate of the
dimension of the center, or better of the homogeneous components of its
positive part, which could allow us to say that the number of independent
central elements cannot be too big: this estimate comes from the
comparison of dim($\z(\u_{\e})\cap\u_{\e,\eta}^+$) with
the multiplicity of
$\e$ in det$H_{\eta}$. 
\prop{\ziz}
Suppose given $\{x_{\a}|\a\in\tPhi_+\}$, with
$x_{\a}\in\u_{\ca,\a}^+$, such that: 

a) $\forall\eta\in Q_+$ $\{x(\ug)|\ug\in\Pp(\eta)\}$ is a basis of
$\u_{\e,\eta}^+$; 

b) $\exists f:J\rightarrow\mz_+$ (with $J\subseteq\tPhi_+$) such that
$x_{\a}^{f({\a})}\in\z(\u_{\e})$ $\forall\a\in J$; 

c) $\forall\eta\in Q_+$ $mult_{\e}{\rm{det}}H_{\eta}\leq\sum_{\a\in J,
m>0}par(\eta-mf(\a)p(\a))$; 

then $\z(\u_{\e})\cap\u_{\e}^+$ is the algebra of
polynomials in $\{x_{\a}^{f(\a)}|\a\in J\}$.
\dim
See \rmcen\ (corollary 3.2.4 and proposition 3.2.5) or equivalently \perf\
(proposition 3.2.12).\finedim
The aim is now to find elements $x_{\a}$'s satisfying the conditions
of proposition \ziz.
\lem{\zelem}
a) $\forall\a\in\tPhi_{+}^{{\rm{re}}}$ the element 
$E_{\a}^{l_{\a}}$ is central in $\u_{\e}$, where $l_{\a}\ugdef{l\over
g.c.d.(l,d_{\a})}$ and $d_{\a}\ugdef d_i$ where $i\in I$ is such that
$\exists w\in W$ with
$\a=w(\a_i)$; 

b) $\forall i\in I_0$, $\forall r>0$ the root vector $E_{(lr\d,i)}$ is
central in $\u_{\e}$. 
\dim 
Part a) is the immediate generalization to the real root
vectors of a classical result by Kac (see \deckac). 
For part b) see the
commutation formulas in
\rcin\ (theorem 5.3.2).\finedim
\rem{\hypjf}
A comparison between proposition \ziz\ and lemma \zelem\ suggests to look
for a set
$J$ containing $J^{\prime}\ugdef\Phi_+^{{\rm{re}}}\cup\{
(lr\d,i)|i\in I_0, r>0\}$ and for a
function $f:J\rightarrow\mz_+$ such that
$f(\a)=l_{\a}$ $\forall\a\in\Phi_+^{{\rm{re}}}$ and
$f((lr\d,i))=1$ $\forall i\in I_0, r>0$; of course following this
suggestion we shall have
$x_{\a}=E_{\a}$ $\forall\a\in J^{\prime}$.
\prop{\multg}
$\forall\eta\in Q_+$ 
the multiplicity of $\e$ in det$H_{\eta}$ is less than or
equal to
\vskip -.4 truecm
$$\sum_{\a\in\Phi_+^{{{\rm{re}}}}\atop m>0}par(\eta-ml_{\a}\a)+
\sum_{r,m>0}mult_{\e}({\rm{det}}\hr)par(\eta-mr\d)$$
\vskip -.4 truecm
where $\hr$ is the matrix defined by 
$\hr\ugdef(E_{(r\d,i)},F_{(r\d,j)})_{ij\in I^r}$ (and $F_{\a}\ugdef
\Omega(E_{\a})$ $\forall\a\in\tPhi_+$).
\dim
For general affine quantum algebras we have that the highest coefficient of
det$H_{\eta}$ is, up to an invertible element of $\mc[q,q^{-1}]$, 
\vskip -.6 truecm
$$\prod_{{\a\in\Phi_+^{{\rm{re}}}\atop m>0}}\left({[m]_{q^{d_{\a}}}\over 
q^{d_{\a}}-q^{-d_{\a}}}\right)^{par(\eta-m\a)}\cdot
\left({\prod_{(r\d,i)\in\tPhi_+^{{\rm{im}}}}
(E_{(r\d,i)},\bar
F_{(r\d,i)})\over\prod_{(r\d,i)\in\tPhi_+^{{\rm{im}}}}A_{ii}^{(r)}}
\right)^{\sum_{m>0}par(\eta-mr\d)}$$
\vskip -.4 truecm
where $\bar F_{(r\d,i)}-A_{ii}^{(r)}F_{(r\d,i)}$ lies in the linear span
of 
$\{F_{(r\d,j)}|j>i\}$ and $(E_{(r\d,i)},\bar F_{(r\d,j)})=0$ $\forall
i>j$. Indeed this is nothing but a reformulation of theorem 2.5.4 of
\rmcen\ (which in this generality does not depend on the
peculiar characteristics of the untwisted algebras, but is valid for all
the affine cases), remarking the connection between the Killing form and
the contravariant form: 
if $x$ and $y$ belong to the linear span of the imaginary root
vectors, $x\in\u_{q,r\d}^+$, $y\in\u_{q,-r\d}^-$, then
$[x,y]=-(x,y)(K_{r\d}-K_{r\d}^{-1})$ (see corollary 7.2.4 of \rcin). 

On the other hand $\forall r>0$ ${\prod_{i\in
I^r} (E_{(r\d,i)},\bar
F_{(r\d,i)})\over\prod_{i\in I^r}A_{ii}^{(r)}}$ is
evidently ${\rm{det}}\hr$ because the matrix of passage from 
$\{\bar F_{(r\d,i)}|i\in I^r\}$ to
$\{F_{(r\d,i)}|i\in I^r\}$ is triangular. The claim then follows remarking
that the multiplicity of
$\e$ in det$H_{\eta}$ is less than or equal to the multiplicity of $\e$ in
the highest coefficient of det$H_{\eta}$.\finedim
\rem{\rcoinc}
Comparing proposition \ziz, remark \hypjf\ and proposition \multg, we see
that for the real part they coincide. So we can
concentrate our attention on the imaginary root vectors, going into the
details of the twisted cases. 
\lem{\twimmult}
In the twisted algebras the multiplicity of $\e$ in
${\rm{det}}\hr$ is given by: 
\vskip -.7 truecm
$$mult_{\e}({\rm{det}}\hr)=\cases{
\#I_0&if $l|r$\cr
1&if $l\not|r$, $2\not|r$, and $l|(2n+1)r$ in case $A_{2n}^{(2)}$ or\cr 
&{\phantom{if $l\not|r$, $2\not|r$, and}} 
$l|\big(\tn-n+1\big)r$ in cases
$A_{2n-1}^{(2)}$ and $E_6^{(2)}$\cr 0&otherwise.}$$ 
\vskip -.3 truecm
\dim
We have that, up to elements which give no
contribution to the multiplicity of $\e$, 
\vskip -.3 truecm
$${\rm{det}}\hr=[r]_q^{\#I^r}\cdot\cases{
[2n+1]_{q^r}&if $2\not|r$ in case $A_{2n}^{(2)}$\cr
[{\tn-n\over k-1}+1]_{q^r}&if $\tk\not|r$ in the other twisted cases\cr
1&otherwise:}$$  
\vskip -.1 truecm
see \rcin\ (theorem 5.3.2, notation 7.1.3, lemma 7.1.4, corollary 7.2.4).
The claim follows immediately remarking that in the case
$D_{\tn}^{(k)}$ we have ${\tn-n\over k-1}=1$.\finedim
\rem{\checerc}
From the comparison between proposition \ziz, remark \hypjf\ and lemma
\twimmult\ we can see that they agree for what concerns the roots
$(r\d,i)$ with $l|r$. But we see also that $J\neq J^{\prime}$: more
precisely, in order to be able to apply proposition \ziz\ we must have
$J=J^{\prime}\cup J^{\prime\prime}$, where $J^{\prime\prime}$ must be of
the form
\vskip -.3 truecm
$$J^{\prime\prime}=\cases{\{(r\d,i_*)|2\not|r,l\not|r,l|(2n+1)r\}&in
case
$A_{2n}^{(2)}$\cr 
\{(r\d,i_*)|2\not|r,l\not|r,l|\big(\tn-n+1\big)r\}&in cases
$A_{2n-1}^{(2)}$ and $E_6^{(2)}$\cr
\emptyset&in cases $D_{n+1}^{(2)}$ and $D_4^{(3)}$;}$$
\vskip -.2 truecm
and $i_*$ must be an element of $I^r$; at the same time we must
have 
$f(\a)=1$ $\forall\a\in J^{\prime\prime}$. This means that we must look
for a nonzero central element in the span of $\{E_{(r\d,i)}|i\in I^r\}$
when $2\not|r,l\not|r,l|(2n+1)r$
in case $A_{2n}^{(2)}$ and 
when $2\not|r,l\not|r,l|(\tn-n+1)r$ in cases $A_{2n}^{(2)}$ and
$E_6^{(2)}$. The results of this section prove that once that we have
found these central elements, we have described
$\z(\u_{\e})\cap\u_{\e}^+$.

{\bf{\zupiu. ${\scriptstyle{\rm{\bf THE\ \ POSITIVE\ \ PART\ \ OF\ \ THE
\ \ CENTER:}}}\ \ \z(\u_{\e})\cap\u_{\e}^+$.}}

The present section is devoted to conclude the program illustrated in
section \ubdz, that is to go into the details in order to exhibit the
central vectors that we are still missing. 
\rem{\qlcnt}
Let $r$ be as described at the end of remark \checerc; then we are looking
for an index $i_*\in I^r$ and for an element $E_{(r\d,i_*)}^*=\sum_{i\in
I^r}A_i^{(r)}E_{(r\d,i)}$ such that: 

a) $E_{(r\d,i_*)}^*$ is linearly independent of $\{E_{(r\d,i)}|i\in
I^r\setminus\{i_*\}\}$ in $\u_{\e}$; 

b) $E_{(r\d,i_*)}^*\in\z(\u_{\e})$.

The above conditions a) and b) can be translated into the following ones: 

a$^{\prime}$) $A_{i_*}^{(r)}\neq 0$ in $\u_{\e}$;

b$^{\prime}$) $(E_{(r\d,i_*)},F_{(r\d,i)})=0$ in $\u_{\e}$ $\forall i\in
I^r$.

Indeed that a) and a$^{\prime}$) are equivalent is obvious; the equivalence
between b) and b$^{\prime}$) is a straightforward consequence of the
commutation relations involving an imaginary root vector (see \rcin,
theorem 5.3.2) and of the already mentioned connection between the bracket
and the Killing form (see the proof of proposition \multg).
But the computations needed to find elements satisfying a$^{\prime}$) and
b$^{\prime}$) have already been carried out in \rcin, so that now we have
just to recall the result.
\prop{\adnd} 
Let $X_{\tn}^{(k)}=A_{2n}^{(2)}$, $2\not|r$, $l\not|r$, $l|(2n+1)r$, and 
let 
\vskip -.3 truecm
$$E_{(r\d,i_*)}^*=[n]_{q^{2r}}E_{(m\d,1)}-\sum_{i\in
I_0\setminus\{1\}}(-1)^r[2]_q[n-i+1]_{q^r}E_{(m\d,i)};$$
\vskip -.2 truecm
then $E_{(r\d,i_*)}^*$ is central in $\u_{\e}$; morover 
$i_*=n$ satisfies the requirement of a$^{\prime}$). 
\dim 
See \rcin, lemma 7.4.1, remark 7.4.2 and propositions 7.4.7 and 7.5.2. As
for the assertion on $i_*$ it is enough to remark that the hypotheses imply
that $l\not|2rn$.
\finedim
\prop{\esd} 
Let $X_{\tn}^{(k)}=A_{2n-1}^{(2)}$ or $E_6^{(2)}$, 
$2\not|r,l\not|r,l|\big(\tn-n+1\big)r$, and let 
\vskip -.3 truecm
$$E_{(r\d,i_*)}^*=\sum_{i\in I^r}(-1)^r[\nu-i+1]_{q^r}\ \ \ {\rm{where}}\
\ \ \nu\ugdef\cases{\tn-n+1=n&in case $A_{2n-1}^{(2)}$\cr 
\tn-n=2&in case
$E_6^{(2)}$;}$$ 
\vskip -.2 truecm
then $E_{(r\d,i_*)}^*$ is central in $\u_{\e}$; morover
$i_*=\nu$ satisfies the requirement of a$^{\prime}$).
\dim
The references given for proposition \adnd\ are still valid, but it is
worth remarking that in these cases the commutation relations involving an
imaginary root vector are exactly those of $A_{\tn-n}^{(1)}$ (see
\rcin, remark 7.4.5), so one can also refer to \rmcen, proposition 3.3.7.
The assertion about $i_*$ is trivial.\finedim
In conclusion the results found above lead to the explicit
description of $\z(\u_{\e})\cap\u_{\e}^+$.
\cor{\cncls}
$\z(\u_{\e})\cap\u_{\e}^+$ is a $\mc$-algebra of
polynomials in an infinite set of variables; more precisely
\vskip -.4 truecm
$$\z(\u_{\e})\cap\u_{\e}^+=\mc[E_{\a}^{l_{\a}}, 
E_{(lr\d,i)},E_{\beta}^*|\a\in\Phi_+^{{\rm{re}}},r>0,i\in I_0,\beta\in
J^{\prime\prime}]$$
\vskip -.2 truecm
where $J^{\prime\prime}$ is as defined in remark \checerc, with $i_*=n$ if
$X_{\tn}^{(k)}=A_{\tn}^{(2)}$, $i_*=2$ in case $E_6^{(2)}$, and the
elements
$E_{\beta}^*$'s are the ones described in propositions \adnd\ and \esd. 
\dim 
This is the natural conclusion of proposition \ziz, lemma \zelem,
proposition \multg, lemma \twimmult\ and propositions \adnd\ and
\esd.\finedim

{\bf{\zue. ${\scriptstyle{\rm{\bf THE\ \ CENTER.}}}$}}

Here we pass from the results obtained till now to the description of the
whole center: ${\z}(\u_{\e})$ will finally turn out to be,
``essentially'', an algebra of polynomials (of course in an infinite number
of variables), with just one relation, regarding its null part. Since we
already know $\z(\u_{\e})\cap\u_{\e}^+$ (then, by
symmetry, we also know $\z(\u_{\e})\cap\u_{\e}^-$), we
are just left with the task of describing
$\z(\u_{\e})\cap\u_{\e}^0$ (which is
trivial) and of understanding how the structure of $\z(\u_{\e})$
can be directly found out from that of its positive, negative and null
parts.
\lem{\zk}
$\z(\u_{\e})\cap\u_{\e}^0$ is the subalgebra of $\u_{\epsilon}^0$
generated by $\{K_i^{l_i},K_{\d}|i\in
I\}$, where
$l_i\ugdef l_{\a_i}$. Remark that while
$\{K_i^{l_i}|i\in I\}$ is a set of algebraically independent elements,
there is a relation between them and
$K_{\d}$: namely, $\prod_{i\in
I}(K_i^{l_i})^{{lr_i\over l_i}}=K_{\d}^l$. 
\dim
The claim is obvious.\finedim
\prop{\zgen}
$\z(\u_{\e})=
(\z(\u_{\e})\cap\u_{\e}^-)\otimes
(\z(\u_{\e})\cap\u_{\e}^0)
\otimes(\z(\u_{\e})\cap\u_{\e}^+)$.
\dim
See \rmcen, theorem 3.4.3, or equivalently \perf, theorem 3.5.6. Remark
that even though the given references relate to the untwisted
affine setting, the proof of the general assertion that we are
dealing with never makes use of the particular form of the untwisted type
algebras but it depends only on the existence of two strings of real root
vectors and on some properties (of the imaginary root vectors and of some
commutation rules) that are common to all the affine situations; the claim
is then valid also in the twisted cases. On the other hand it is worth
noticing that this result is not true for the quantum algebras of finite
type (these are the quantization of the enveloping algebras of simple
finite dimensional complex Lie algebras), where there are also the Casimir
elements, that is central elements which cannot be decomposed as algebraic
combinations of ``positive'', ``negative'' and ``null'' central
elements.\finedim
\teo{\centrofine}
Let $\u_{\e}$ be the specialization at $\e$ of an affine
quantum algebra of twisted type $X_{\tn}^{(k)}$, with $\e\in\mc$
primitive
$l^{{\rm{th}}}$ of 1 and $l$ odd integer bigger than $k$. Then the center
of $\u_{\e}$ is  
\vskip -.7 truecm
$$\z(\u_{\e})\!\!=\!\![E_{\a}^{l_{\a}}, 
E_{(lr\d,i)},E_{\beta}^*,F_{\a}^{l_{\a}}, 
F_{(lr\d,i)},F_{\beta}^*,K_j^{l_j},K_{\d}|\a\!\in\Phi_+^{{\rm{re}}},r>0,i\in
I_0,j\in I,\beta\in J^{\prime\prime}]/(P_Z)$$
\vskip -.2 truecm
where $J^{\prime\prime}$ and $E_{\beta}^*$ are those of corollary
\cncls, $F_{\beta}^*\ugdef\Omega(E_{\beta}^*)$ and 
$P_Z\ugdef K_{\d}^l-\prod_{i\in I}(K_i^{l_i})^{{lr_i\over l_i}}$. 
\dim
The theorem is the straightforward consequence of corollary \cncls, lemma
\zk\ and proposition \zgen.\finedim

{\bf{\ref. References.}}

\beck\ Beck, J., {\it{Convex bases of PBW type for quantum affine
algebras}},  Commun. Math. Phys. {\bf 165} (1994), 193-199. 

\rmcen\ Damiani, I., {\it{The 
highest coefficient of det$H_{\eta}$ and the center of the specialization 
at odd roots of unity for untwisted affine quantum algebras}}, 
J. Algebra {\bf 186} (1996), 736-780. 

\rcin\ Damiani, I., {\it{The $R$-matrix for (twisted) affine quantum
algebras}}, Proceedings of the International Conference on Representatin
Theory, June 29-July 3, 1998, East China Normal University, Shanghai,
China, China Higher Education Press \&\ Springer-Verlag, Beijing (2000),
89-144. 

\perf\ Damiani, I., {\it{Untwisted affine quantum algebras: the highest
coefficient of det$H_{\eta}$ and the center at odd roots of 1}}, tesi di
perfezionamento, Scuola Normale Superiore - Pisa (1996). 

\deckac\ De Concini, C., Kac, V.G., {\it{Representations of quantum
groups at roots of 1}}, Progr. in Math. {\bf 92}, Birk\"auser (1990),
471-506

\KM\ Kac, V.G., {\it{Infinite Dimensional Lie Algebras}}, 
Birkh\"auser Boston, Inc., 
USA (1983). 

\lusg\ Lusztig, G., {\it{Quantum deformations of certain simple modules
over enveloping algebras}}, Adv. in Math. {\bf 70} (1988), 237-249. 

\lusz\ Lusztig, G., {\it{Quantum groups at roots of 1}}, Geom. Ded. 
{\bf 35} (1990), 89-113. 

\tan\ Tanisaki, T., {\it{Killing forms, Harish-Chandra isomorphisms
and universal R-matrices for quantum algebras}}, in Infinite
Analysis Part B, Adv. Series in Math. Phys., vol 16, 1992, p. 941-962.

\bye